\documentclass[a4paper,12pt,leqno]{amsart}
\usepackage{amsmath,amssymb,mathrsfs}
\usepackage[matrix,arrow]{xy} 
\usepackage{xcolor,hyperref}
\hypersetup{colorlinks=true,citecolor=black,linkcolor=black}
\newcommand{\h}{\hbox}
\newcommand{\q}{\quad}
\newcommand{\nin}{\noindent}

\newcommand{\ms}{\par\medskip}
\newcommand{\sk}{\par\smallskip}

\newcommand{\msn}{\par\medskip\noindent}

\newcommand{\ges}{\geqslant}
\newcommand{\les}{\leqslant}
\newcommand{\one}{\hskip1pt}

\newcommand{\mcup}{\hbox{$\bigcup$}}
\newcommand{\msum}{\hbox{$\sum$}}

\newcommand{\G}{{\mathcal G}}
\newcommand{\K}{{\mathcal K}}
\newcommand{\Rc}{{\mathcal R}}

\newcommand{\Q}{{\mathbb Q}}
\newcommand{\C}{{\mathbb C}}
\newcommand{\N}{{\mathbb N}}
\newcommand{\R}{{\mathbb R}}

\newcommand{\Z}{{\mathbb Z}}
\newcommand{\ob}{{\mathbf 1}}

\newcommand{\Gr}{{\rm Gr}}

\newcommand{\aTj}{\alpha^{\rm\hskip-1pt Tj}}
\newcommand{\Sb}{{\bf S}}
\newcommand{\dfTj}{\delta_{\hskip-1pt f}^{\rm Tj}}
\newcommand{\Cf}{C_{\hskip-1pt f}}
\newcommand{\If}{I_{\hskip-1pt f}}
\newcommand{\Tf}{T_{\!f}}
\newcommand{\Gf}{\G_{\hskip-1pt f}}
\newcommand{\Hsf}{H''_{\hskip-1pt f}}
\newcommand{\dti}{\dd_t^{-1}}
\newcommand{\dt}{\dd_t}
\newcommand{\al}{\alpha}
\newcommand{\be}{\beta}

\newcommand{\la}{\lambda}
\newcommand{\La}{\Lambda}

\newcommand{\om}{\omega}
\newcommand{\Om}{\Omega}
\newcommand{\dd}{\partial}
\newcommand{\ddd}{{\rm d}}
\newcommand{\Sp}{{\rm Sp}}
\newcommand{\SpT}{{\rm Sp}^{\mhs\rm Tj}}
\newcommand{\Spf}{{\rm Sp}_{\mhs f}}
\newcommand{\SpTf}{{\rm Sp}\mhs^{\rm Tj}_f}
\newcommand{\Varf}{{\rm Var}\mhs_f}
\newcommand{\VarTf}{{\rm Var}\mhs^{\rm Tj}_f}
\newcommand{\avTf}{{\rm av}\mhs^{\rm Tj}_{\mhs f}}

\newcommand{\eq}{\,{=}\,}
\newcommand{\defs}{\,{:=}\,}
\newcommand{\nes}{\,{\ne}\,}
\newcommand{\ins}{\,{\in}\,}
\newcommand{\sst}{\,{\subset}\,}
\newcommand{\stm}{\,{\setminus}\,}
\newcommand{\gess}{\,{\ges}\,}
\newcommand{\less}{\,{\les}\,}
\newcommand{\sgt}{\,{>}\,}
\newcommand{\slt}{\,{<}\,}

\newcommand{\pl}{\one {+}\one}
\newcommand{\mi}{\one {-}\one}
\newcommand{\bl}{\bigl}
\newcommand{\br}{\bigr}

\newcommand{\ssb}{\raise.15ex\h{${\scriptscriptstyle\bullet}$}}
\newcommand{\ssc}{\,\raise.15ex\h{${\scriptstyle\circ}$}\,}
\newcommand{\onto}{\twoheadrightarrow}

\newcommand{\simto}{\,\,\rlap{\hskip1.5mm\raise1.4mm\hbox{$\sim$}}\hbox{$\longrightarrow$}\,\,}
\newcommand{\mhs}{\hskip-1pt}

\newcommand{\pn}{\par\noindent}

\makeatletter
\renewcommand\section{\@startsection{section}{1}{0pt}{-2ex plus -1ex minus -.2ex}{-2.3ex plus-.2ex}{\normalfont\bfseries}}
\makeatother
\theoremstyle{plain}
\newtheorem{thm}{Theorem}[section]

\newtheorem{prop}{Proposition}[section]

\theoremstyle{definition}
\newtheorem{rem}{Remark}[section]

\newtheorem{iconj}{Conjecture}

\begin{document}
\title[Generalized Hertling conjecture]
{Some remarks on the generalized Hertling conjecture for Tjurina spectrum}
\author[S.-J. Jung]{Seung-Jo Jung}
\address{S.-J. Jung : Department of Mathematics Education, and Institute of Pure and Applied Mathematics, Jeonbuk National University, Jeonju, 54896, Korea}
\email{seungjo@jbnu.ac.kr}
\author[I.-K. Kim]{In-Kyun Kim}
\address{I.-K. Kim : June E Huh Center for Mathematical Challenges, Korea Institute for Advanced Study, 85 Hoegiro Dongdaemun-gu, Seoul 02455, Korea}
\email{soulcraw@kias.re.kr}
\author[M. Saito]{Morihiko Saito}
\address{M. Saito : RIMS Kyoto University, Kyoto 606-8502 Japan}
\email{msaito@kurims.kyoto-u.ac.jp}
\author[Y. Yoon]{Youngho Yoon}
\address{Y. Yoon : Department of Mathematics, Chungbuk National University, Cheongju-si, Chungcheongbuk-do, 28644, Korea}
\email{mathyyoon@gmail.com}
\thanks{This work was partially supported by National Research Foundation of Korea (the first author: NRF-2021R1C1C1004097, the second author: NRF-2023R1A2C1003390 and NRF-2022M3C1C8094326, and the fourth author: RS-2023-00245670).}
\begin{abstract} We study the original version of the generalized Hertling conjecture on the variance of the Tjurina spectral numbers, which was proposed by Shi, Wang, and Zuo, and provide a sufficient condition for the original conjecture to fail, employing a theorem of Hertling in an essential way. We calculate certain examples using some codes in Singular.
\end{abstract}
\maketitle

\part*{Introduction} \label{intr}
\nin
Let $f\ins\C\{x\}$ be a convergent power series having an isolated singularity with $f(0)\eq0$, where $x\eq(x_1,\dots,x_n)$ is the coordinate system of $(\C^n,0)$. One can define the {\it Steenbrink spectrum\one} $\Spf(t)\eq\msum_{i=1}^{\mu}\,t^{\one\al_i}$ as the Poincar\'e polynomial of the $V$-filtration on the Jacobian ring $\C\{x\}/(\dd f)$, that is,
\begin{equation} \label{1}
|\If^{\al}|=\dim_{\C}\Gr^{\al}_V\bl(\C\{x\}/(\dd f)\br)\q\h{with}\q\If^{\al}\eq\{i\in\If\mid\al_i\eq\al\}.
\end{equation}
Here $(\dd f)\sst\C\{x\}$ is the Jacobian ideal, $V$ is the quotient filtration of the $V$-filtration on the Brieskorn lattice (or the {\it microlocal\one} $V$-filtration on $\C\{x\}$) indexed by $\Q$, and $\If\defs\{1,\dots,\mu\}$ with $\mu$ the Milnor number, see \cite{SS}, \cite{Va1}, and also \cite{bl}, \cite{JKSY0}, \cite{JKSY2}. We may assume that $f$ is a polynomial by the finite determinacy, as is well known. The spectral numbers $\al_i$ are assumed to be {\it weakly increasing.}
\sk
The {\it Tjurina subspectrum\one} (or {\it spectrum\one}) $\SpTf(t)=\msum_{j=1}^{\tau}\,t^{\one\aTj_j}=\msum_{i\in\Tf}\,t^{\one\al_i}$ with $\Tf\sst\If$ is defined by
\begin{equation} \label{2}
|\Tf^{\al}|=\dim_{\C}\,\Gr^{\al}_V\bl(\C\{x\}/(\dd f,f)\br)\q\h{with}\q\Tf^{\al}\eq\If^{\al}\cap\Tf\eq\If^{\al}\cap [1,i_{\al}],
\end{equation}
with $i_{\al}\ins\Tf^{\al}$, where the $\aTj_j$ are weakly increasing, see \cite{JKY}, \cite{JKSY2}. Note that $|\Tf|\eq\SpTf(1)\eq\tau$ with $\tau$ the Tjurina number of $f$. Set $\Cf\defs\If\stm\Tf$, and
\begin{equation} \label{3}
\Spf^C(t)\defs\Spf(t)\mi\SpTf\eq\msum_{i\in\Cf}\,t^{\one\al_i}\eq\msum_{k=1}^{\mu-\tau}\,t^{\one\al^C_k},
\end{equation}
where the $\al^C_k$ are assumed to be weakly increasing. By definition the fractional polynomial $\Spf^C(t)$ is the Poincar\'e polynomial of the {\it image\one} of the endomorphism
$$[f]\ins{\rm End}_{\C\{x\}}\bl(\C\{x\}/(\dd f)\br),$$
\pn with $[f]\in\C\{x\}/(\dd f)$ the class of $f$, and is called the {\it complemental part\one} to the Tjurina spectrum. We call the numbers $\al^C_k$ ($k\ins[1,\mu{-}\tau]$) and $\al_i$ ($i\ins\Cf$) the {\it missing spectral numbers\one} of $f$. Some of these missing spectral numbers may also belong to the Tjurina spectrum. It has been proved that the missing spectral numbers have a {\it graded symmetry,} see \cite{JKSY2}. Note also that the Tjurina spectrum is {\it unstable under a $\one\tau$-constant deformation,} see \cite[Remark 1]{JKSY2}.
\sk
It has been conjectured that the {\it variance\one} of the Tjurina spectral numbers is bounded by the width of the Tjurina spectrum (that is, the difference between the maximal and minimal Tjurina spectral numbers) divided by 12, as conjectured by Shi, Wang, and Zuo \cite{SWZ25}, \cite{SWZ26} (see also Section\,\,\ref{S3} below). In this paper we show that this original version of the conjecture does not necessarily hold. They now conjecture the inequality only for the case of $\tau$-max form, excluding the semi-weighted-homogeneous case, see \cite[Conjecture 1.3]{SWZ26}. It seems then rather interesting to consider the difference $\dfTj$ of the two numbers (see \eqref{3.2} below) as a {\it subtle analytic invariant\one} of a non-weighted-homogeneous polynomial with an isolated singularity at 0, {\it measuring the complexity of singularity in some sense.}
\sk
The simplest example (with respect to $\mu$) such that $\dfTj\sgt0$ seems to be $f\eq x^7\pl y^7\pl x^5y^5$ with $\mu\eq36$, $\mu{-}\tau\eq1$, and modality 10 (combining \cite{Ga} and \cite{Va3}). It is interesting that any $\mu$-constant deformation of $x^a\pl y^b$ with $b\less a$, $b\less 7$, and $\mu\nes\tau$ seems to have non-positive $\dfTj$ except the case $a\eq b\eq 7$ with $\mu{-}\tau\eq1$ as far as examples are computed, see Section\,\,\ref{S6} below. One may make the following.

\begin{iconj} \label{Con1}
The original version of the generalized Hertling conjecture holds as long as $f$ {\it is not semi-weighted-homogeneous.}
\end{iconj}

Counterexamples to the original conjecture are restricted to semi-weighted-homogeneous polynomials with small $(\mu{-}\tau)/\mu$, and the latter number for non-semi-weighted-homogeneous ones does not seem to be as small as that for the counterexamples in general.

\sk
In Part~1, after reviewing some basics of spectrum, we provide a sufficient condition for the generalized Hertling conjecture on the variance of Tjurina spectrum to fail. In Part~2 we confirm the conclusions of some conceptual theorems by calculating explicit examples.

\tableofcontents
\numberwithin{equation}{section}

\part{Generalized Hertling conjecture on variance} \label{Pa1}
\nin
In this part, after reviewing some basics of spectrum, we provide a sufficient condition for the generalized Hertling conjecture on the variance of Tjurina spectrum to fail.

\section{Spectrum of hypersurface isolated singularities} \label{S1}
Let $f\ins\C\{x\}$ be a convergent power series of $n$ variables having an isolated singularity at 0, where $f(0)\eq0$. There is a canonical mixed Hodge structure on the vanishing cohomology $H^{n-1}(F_{\!f},\Q)$ (using the natural coordinate $t$ of $\C$) with $F_{\!f}$ the Milnor fiber of $f$, see \cite{St}. (This can be defined by using the vanishing cycle functor $\varphi_f$ of mixed Hodge modules \cite{mhm}.) The {\it spectrum\one} $\Sp_f(t)\eq\msum_{i\in\If}\,t^{\al_i}$ of $f$ (where $\If\defs\{1,\dots,\mu\}$ with $\mu$ the Milnor number of $f$) is defined by
\begin{equation} \label{1.1}
\aligned&\#\{i\ins\If\mid\al_i\eq\al\}=\dim_{\C}\Gr_F^pH^{n-1}(F_{\!f},\C)_{\la}\\&\q\q\q\q\q\h{with}\q p\defs[n{-}\al],\,\la\defs e^{-2\pi \sqrt{-1}\al}.\endaligned
\end{equation}
Here $F$ is the Hodge filtration of the mixed Hodge structure, and $_{\la}$ denotes the $\la$-eigenspace of the action of the semi-simple part $T_s$ of the Jordan decomposition of the monodromy $T$ (which is the inverse of the Milnor monodromy), see \cite{St} and also \cite{DS}. The rational numbers $\al_1,\dots,\al_{\mu}$ are contained in $(0,n)$, and are called the {\it spectral numbers\one} of $f$. These are indexed weakly increasingly.
We have the {\it symmetry\one} of spectral numbers (see \cite{St}):
\begin{equation} \label{1.2}
\al_i\pl\al_j\eq n\q\h{if}\,\,\,i{+}j\eq\mu{+}1.
\end{equation}
Here we use the well-known assertion that the weight filtration is given by the monodromy filtration, shifted by $n$ or $n{-}1$ depending on whether the monodromy eigenvalue is 1 or not. Note also that
\begin{equation} \label{1.3}
\h{The spectrum stays invariant under $\mu$-constant deformations.}
\end{equation}
This follows from \cite{Va2} (or \cite{DMST}).
\sk
Let $\Hsf$ be the Brieskorn lattice. This is contained in the Gauss-Manin system $\Gf$, which has the $V$-filtration of Kashiwara and Malgrange indexed by $\Q$ (which were originally indexed by $\Z$), see for instance \cite{bl}. Set $F_p\Gf\defs\dd_t^p\Hsf\sst\Gf$ for $p\ins\Z$. It is well known that there are isomorphisms
\begin{equation} \label{1.4}
\aligned&F^p\!H^{n-1}(F_{\!f},\C)_{\la}=F\!_{-p-j}\Gr_V^{\al+j}\Gf\\&\q\h{for any}\,\,\,\,p,j\ins\Z,\,\al\ins(0,1],\,\la\eq e^{-2\pi\sqrt{-1}\al},\endaligned
\end{equation}
see \cite{SS}, \cite{Va1} (and also \cite[\S3.4]{gm1}). Moreover, in the notation of the Introduction we have the isomorphism
\begin{equation} \label{1.5}
\Gr_F^0\Gf=\Hsf/\dd_t^{-1}\Hsf=\C\{x\}/(\dd f),
\end{equation}
using the coordinates $x_1,\dots,x_n$ which give a trivialization of $\Om_{\C^n}^n$ by $\ddd x_1{\wedge}\cdots{\wedge}\ddd x_n$. So the spectrum can be defined also as in \eqref{1}.

\section{Hertling conjecture on variance} \label{S2}
The {\it variance\one} of the spectral numbers $\al_1,\dots,\al_{\mu}$ is defined by
\begin{equation} \label{2.1}
\Varf:=\tfrac{1}{\mu}\,\msum_{i\in\If}\bl(\al_i\mi\tfrac{n}{2}\br)^2,
\end{equation}
since ${\rm av}\!_f\defs\tfrac{1}{\mu}\msum_{i\in\If}\al_i\eq\tfrac{n}{2}$ using the symmetry of spectral numbers $\al_i$, see \eqref{1.2}. Hertling's conjecture asserts the inequality
\begin{equation} \label{2.2}
\Varf\les\tfrac{1}{12}(\al_{\mu}\mi\al_1),
\end{equation}
see \cite{He}. This has been shown in the curve case, see \cite{irr} for the irreducible case, \cite{Bre1} for the non-degenerate case, and \cite{Bre2} for the general curve case. (The last paper does not seem to be published yet.) Note that the equality holds in the weighted homogeneous case, see \cite{He}.

\section{Generalized Hertling conjecture} \label{S3}
In the notation of the introduction of \cite{JKSY2}, the variance of the Tjurina spectral numbers $\al_i$ ($i\ins\Tf$) is defined by
\begin{equation} \label{3.1}
\VarTf:=\tfrac{1}{\tau}\,\msum_{i\in\Tf}\bl(\al_i\mi\avTf\br)^2\q\h{with}\q\avTf:=\tfrac{1}{\tau}\,\msum_{j\in\Tf}\,\al_j.
\end{equation}
The original version of the generalized Hertling conjecture on the variance of the Tjurina spectrum \cite{SWZ25}, \cite{SWZ26} claims the inequality
\begin{equation} \label{3.2}
\dfTj:=\VarTf\mi\tfrac{1}{12}\bl(\aTj_{\tau}\mi\aTj_1\br)\les 0,
\end{equation}
where $\aTj_{\tau}$ is the maximal Tjurina spectral number and $\aTj_1\eq\al_1$, see \cite[Theorem 2]{JKSY2}.
\sk
Applying Hertling's theorem mentioned at the end of Section~\ref{S2}, we prove the theorem below which implies that the conjecture fails in certain cases.

\begin{thm} \label{T3.1}
Assume $f$ is semi-weighted-homogeneous with $\mu\nes\tau$ and either $\al_{\mu}\mi\al_1\less 2$ or more generally $\avTf\less{\rm av}\!_f$. Then the inequality {\rm\eqref{3.2}} does not hold if
\begin{equation} \label{3.3}
\tfrac{\mu}{12}(\al_{\mu}\mi\aTj_{\tau})\ges(\mu{-}\tau)\al_{\mu}^2.
\end{equation}
\end{thm}

\begin{proof}
Note first that the inequality $\avTf\less{\rm av}\!_f$ follows from the condition $\al_{\mu}\mi\al_1\less 2$ using the inequality in \cite[(3)]{JKSY2}. We have the inequalities
\begin{equation} \label{3.4}
\aligned\msum_{i\in\Tf}\bl(\al_i\mi\tfrac{1}{\tau}\one\msum_{j\in\Tf}\,\al_j\br)^2&=\msum_{i\in\Tf}\,\al_i^2-\tau\bl(\one\tfrac{1}{\tau}\one\msum_{j\in\Tf}\,\al_j\br)^2\\&>\msum_{i\in\Tf}\,\al_i^2-\mu\bl(\one\tfrac{1}{\mu}\one\msum_{j\in\If}\,\al_j\br)^2\\&\ges\tfrac{\mu}{12}(\al_{\mu}\mi\al_1)-(\mu{-}\tau)\al_{\mu}^2\\&>\tfrac{\tau}{12}(\aTj_{\tau}\mi\al_1)+\tfrac{\mu}{12}(\al_{\mu}\mi\aTj_{\tau})-(\mu{-}\tau)\al_{\mu}^2.\endaligned
\end{equation}
The second inequality follows from the assertion in \cite{He} mentioned at the end of Section\,\,\ref{S2}, which implies that the equality holds in \eqref{2.2} for semi-weighted-homogeneous polynomials using \eqref{1.3}. So the assertion follows. This finishes the proof of Theorem\,\,\ref{T3.1}.
\end{proof}

\begin{rem} \label{R3.1}
The condition $\avTf\less{\rm av}\!_f$ seems to be always satisfied as far as examples are calculated.
\end{rem}

\begin{rem} \label{R3.2}
Let $f'$ be a $($virtual$\one)$ $\mu$-constant deformation of $f$ such that its Tjurina number $\tau'$ coincides with $\tau$ and $\SpTf\mi\SpT_{\mhs f'}\eq t^{\be}\mi t^{\be'}$. Let $\al_{\tau}^{\prime\,{\rm Tj}}$ be the maximal Tjurina spectral number of $f'$. In the case $\be\sgt\be'$ and $\aTj_{\tau}\mi\al_{\tau}^{\prime\,{\rm Tj}}\eq\be\mi\be'$, we have
\begin{equation} \label{3.5}
\dfTj<\delta_{\mhs f'}^{\rm Tj}\q\h{if}\q\be\pl\be'\les\avTf\pl{\rm av}\mhs^{\rm Tj}_{\mhs f'}\pl\tfrac{\tau}{12}.
\end{equation}
On the other hand, if $\be\sgt\be'$ and $\aTj_{\tau}\eq\al_{\tau}^{\prime\,{\rm Tj}}$, we have
\begin{equation} \label{3.6}
\dfTj>\delta_{\mhs f'}^{\rm Tj}\q\h{if}\q\be\pl\be'\ges\avTf\pl{\rm av}\mhs^{\rm Tj}_{\mhs f'}.
\end{equation}
These follow from the first equality of \eqref{3.4}, since
\begin{equation} \label{3.7}
\aligned\msum_{i\in\Tf}\,\al_i^2\mi\msum_{i\in T_{\!f'}}\,\al_i^2&=(\be{-}\be')(\be{+}\be'),\\ \tau(\avTf)^2\mi\tau({\rm av}\mhs^{\rm Tj}_{\mhs f'})^2&=(\be{-}\be')(\avTf\pl{\rm av}\mhs^{\rm Tj}_{\mhs f'}).\endaligned
\end{equation}
\sk
By \eqref{3.6} the proof of the inequality \eqref{3.2} is {\it reduced\one} to the case where $\Tf\cap\If^{\ges\al_1+1}$ is {\it consecutive\one} (that is, equal to $\Z\cap[a,b]$ for some $a,b\in\Z$) assuming $\al_{\mu}{-}\al_1\slt2$.
\end{rem}

\section{A slight extension} \label{S4}
Sometimes it is useful to consider the inequality \eqref{3.2} for a subset $T\sst\If$ which does not necessarily coincide with $T_g$ for a $\mu$-constant deformation $g$ of $f$. The inequality \eqref{3.2} can be generalized to
\begin{equation} \label{4.1}
\delta_T:={\rm Var}_T-\tfrac{1}{12}(\al_T^{\rm max}\mi\al_T^{\rm min})\les0.
\end{equation}
Here ${\rm Var}_T\defs\tfrac{1}{\tau}\msum_{i\in T}\one (\al_i\mi{\rm av}_T)^2$ with ${\rm av}_T\defs\tfrac{1}{\tau}\one\msum_{i\in T}\one\al_i$, $\tau\defs|T|$, and $\al_T^{\rm max}\defs\max\{\al_i\,|\,i\ins T\}$, etc. We have the following.

\begin{prop} \label{P4.1}
Let $T'\defs T\stm\{i_0\}$ with $i_0\ins T$ such that $\al_{T'}^{\rm max}\eq\al_T^{\rm max}$, $\al_{T'}^{\rm min}\eq\al_T^{\rm min}$. Then the inequality {\rm\eqref{4.1}} holds for $T'$ if it holds for $T$ and
\begin{equation} \label{4.2}
(\al_{i_0}\mi{\rm av}_T)^2\gess\tfrac{1}{12}(\al_T^{\rm max}\mi\al_T^{\rm min}).
\end{equation}
\end{prop}

\begin{proof}
We have
\begin{equation} \label{4.3}
\msum_{i\in T}\,\bl(\al_i\mi\tfrac{1}{\tau}\msum_{i\in T}\,\al_i\br)^2=\msum_{i\in T}\,\al_i^2\mi\tfrac{1}{\tau}\bl(\msum_{i\in T}\,\al_i\br)^2,
\end{equation}
and similarly with $T$ replaced by $T'$ and $\tau$ by $\tau'\defs|T'|\eq\tau\mi1$. Here we may assume ${\rm av}_T\eq0$ replacing the $\al_i$ with $\al_i\mi{\rm av}_T$. By the assumption \eqref{4.2} we then get
\begin{equation} \label{4.4}
\tau\delta_T\gess\tau'\delta_{T'},
\end{equation}
since $\bl(\msum_{i\in T'}\,\al_i\br)^2\gess0$. So Proposition\,\,\ref{P4.1} follows.
\end{proof}

\begin{rem} \label{R4.1}
It has been conjectured by K.\,Saito that the difference between the maximal and minimal exponents of a hypersurface isolated singularity is at most 1 if and only if the singularity is rational double or simple elliptic or cusp, see for instance \cite{SWZ25}, \cite{SWZ26}. This conjecture, however, does not hold with their formulation in general. The example just below corresponds to choosing a non-standard order of $\R$ (see \cite[\S3.9]{bl}), and is simpler than the one in \cite[\S4.4]{bl}. (They say that such a non-standard order of $\R$ is never used, but this denial never follows from their formulation.)
\sk
Assume for instance $f\eq x^5\pl y^4\pl x^3y^2$. There are free generators $\eta_i$ of $\Gf$ over $\K$ such that $\dt t\eta_i\eq\al_i\eta_i$ and $\Sb(\eta_i,\eta_j)\eq\delta_{i,\one\mu+1-j}\one\dt^{-n}$ ($i,j\ins[1,\mu]$) with $\mu\eq12$ in the notation of \cite[\S1.2]{JKSY2}, since $\dim H^1(F_{\!f},\C)_{\la}\less 1$. Note that the spectral numbers of $f$ are given by
\begin{equation} \label{4.5}
\bl\{\tfrac{4p+5q}{20}\mid(p,q)\ins[1,4]{\times}[1,3]\br\}.
\end{equation}
The Brieskorn lattice $\Hsf$ is generated over $\Rc$ by the $\om_i$ ($i\ins[1,\mu]$) with
\begin{equation} \label{4.6}
\om_1\defs\eta_1\pl\dt\eta_{\mu},\q\om_i\defs\eta_i\,\,\,(i\ins[2,\mu]),
\end{equation}
using the theory of opposite filtrations \cite{bl}. Here $\eta_1,\eta_{\mu}$ are replaced by $c\one\eta_1,c^{-1}\eta_{\mu}$ for some $c\ins\C^*$. This gives a {\it very good section,} since $\Sb(\om_i,\om_j)\eq\delta_{i,\one\mu+1-j}\one\dt^{-n}$, see \cite{bl}, \cite{bl2}. One can however replace $\om_{\mu}\defs\eta_{\mu}$ with
\begin{equation} \label{4.7}
\om'_{\mu}\defs{-}\dti\eta_1\eq\eta_{\mu}\mi\dti\om_1.
\end{equation}
Its image in $\Om_f$ is unchanged. Setting $\om'_i\defs\om_i$ ($i\ins[1,\mu{-}1]$), one has $\Sb(\om'_i,\om'_j)\eq\delta_{i,\one\mu+1-j}\one\dt^{-n}$ up to sign. This gives also a good section, but it is {\it not\one} very good. Indeed, the exponents of the modified section for the $\Rc$-submodule $\Rc\om'_1\pl\Rc\om'_{\mu}$ are $\al_1{+}1\eq\tfrac{29}{20}$ and $\al_{\mu}{-}1\eq\tfrac{11}{20}$, since
\begin{equation} \label{4.8}
\dt t\om'_{\mu}\eq\tfrac{29}{20}\om'_{\mu},\q\dt t\om'_1\eq\tfrac{11}{20}\om'_1\pl c'\eta_1\eq\tfrac{11}{20}\om'_1\mi c'\dt\om'_{\mu},
\end{equation}
for some $c'\ins\C^*$. The minimal and maximal exponents are then $\tfrac{11}{20}$ and $\tfrac{29}{20}$ in view of \eqref{4.5}.
\end{rem}

\part{Explicit calculations} \label{Pa2}
\nin
In this part we confirm the conclusions of some conceptual theorems by calculating explicit examples.

\section{Ordinary {\it m}-ple point case} \label{S5}
Assume $f$ defines an ordinary $m$-ple point, that is, the leading term $f_1$ of $f\eq\msum_{\be\ges1}f_{\be}$ is a homogeneous polynomial of degree $m$ having an isolated singularity. We have $\mu\eq(m{-}1)^n$ and $\al_{\mu}\mi\aTj_{\tau}\gess\tfrac{1}{m}$. By Theorem\,\,\ref{T3.1} the inequality \eqref{3.2} then fails if
\begin{equation} \label{5.1}
(m{-}1)^n\gess12\one m\one n^2(\mu{-}\tau),
\end{equation}
assuming $\avTf\less{\rm av}\!_f$. In the case $\mu{-}\tau$ is fixed (for instance if $\mu{-}\tau\eq1$), these conditions are satisfied when $m$ is sufficiently large. In the case $n\eq2$ with $\mu{-}\tau\eq1$, the inequality \eqref{3.2} fails actually for $m\gess 7$ according to a computation in Section\,\,\ref{S6} just below.

\section{Semi-weighted-homogeneous case} \label{S6}
The situation is slightly different in the general semi-weighted-homogeneous case. Using Singular (see \cite{Sing}), one can examine the inequality \eqref{3.2} for $n\eq2$ as follows.
\ms
\vbox{\footnotesize\sf\verb#LIB"sing.lib"; ring R=0, (x,y), ds; int a,b,c,d; a=7; b=7; c=1; d=1;#
\sk
\verb#int p, q, i, j, t, N; poly f, Sm, Vt, av, sm, ai, bi, u, Max, X;#
\sk
\verb#u=1; ai=u/a; bi=u/b; p=a-1-c; q=b-1-d; f=x^a+y^b+x^p*y^q;#
\sk
\verb#sm=(a-1)*(b-1)-((2*a-1-c)*ai+(2*b-1-d)*bi)/2*c*d;#
\sk
\verb#Sm=0; Vt=0; N=0; Max=0; t=(a-1)*(b-1)-c*d; av=sm/t;#
\sk
\verb#for(i=1; i<=a-1; i++) {for(j=1; j<=b-1; j++) {if(i<a-c||j<b-d) {#
\sk
\verb#X=i*ai+j*bi; Sm=Sm+X; Vt=Vt+(X-av)^2; N=N+1; if(X>Max) {Max=X;}}}}#
\sk
\verb#if(N!=t || Sm/t!=av || tjurina(f)!=t) {printf("Error!");}#
\sk
\verb#Vt/t - (Max-ai-bi)/12;#}
\msn
One can change the definitions of the positive integers $a,b,c,d$ as long as they satisfy the conditions:
$$c\slt\tfrac{a}{2},\q d\slt\tfrac{b}{2},\q \tfrac{a-1-c}{a}\pl\tfrac{b-1-d}{b}\sgt1.$$
\par\nin These conditions should imply the equality $\mu{-}\tau\eq cd$. Here it is {\it insufficient\one} to assume only the last inequality. The inequality \eqref{3.2} seems to hold, that is, the last output is a non-positive number, if $b\less a$ and $b\less 7$, except the case $a\eq b\eq 7$ with $c\eq d\eq1$. It may hold for any $a,b$ if $c,d$ are sufficiently large.

\section{Some technical difficulties} \label{S7}
In \cite{SWZ26} the computer program Singular is used for the proof of the inequality \eqref{3.2} in the case the modality ${\rm mod}_f$ is at most 3. Here one has to take a rational point on each connected stratum of a stratification of the base space of the miniversal $\mu$-constant deformation of $f$ such that the Tjurina spectrum is constant on each stratum. This seems very difficult, since the Tjurina spectrum is not stable by $\tau$-constant deformations. It may be simpler to try to prove the inequality \eqref{3.2} for all the {\it possible candidates\one} for the Tjurina spectrum for each {\it possible\one} $\tau$, and show that it does not occur in the case the inequality does not hold. However, this method may have some difficulty if $|\If^{\ges\al_1+1}|$ is quite big (see the inclusion in (3) of \cite[Theorem~2]{JKSY2}), since one may get many possibilities of missing spectral numbers. In this case the last part of Remark\,\,\ref{R3.2} is sometimes helpful, and in the case where the last assumption of \cite[Theorem~2]{JKSY2} is satisfied, one can consider for instance the following.
\ms
\vbox{\footnotesize\sf\verb#LIB"sing.lib"; LIB"gmssing.lib"; ring R=0,(x,y),ds; poly f,Av,Sm,Vt; #
\sk
\verb#list sp; int i,j,k,m,t,s,c,A,r,mu; f=x^7+y^7; A=10; sp=spectrum(f);#
\sk
\verb#s=size(sp[2]); mu=milnor(f); matrix S[1][mu]; m=1; for(i=1; i<=s; i++)#
\sk
\verb#{for(j=1; j<=sp[2][i]; j++) {S[1,m]=sp[1][i]+1; m++;}} m--; if(m!=mu)#
\sk
\verb#{sprintf(" m Error"); exit;} t=tjurina(f); for (k=m-1; S[1,k]>S[1,1]+1;#
\sk
\verb#k--) {;} k++; if (k>t-A+1){A=t-k+1; sprintf(" A replaced by %s !",A);}#
\sk
\verb#sprintf(" mu_f=%s, tau_f=%s, 1+al_1=%s, al_%s=%s, al_%s=%s",m,t,S[1,1]+1,#
\sk
\verb#k-1,S[1,k-1],k,S[1,k]); for (r=0; r<=A; r++) {for (j=m-t; j>0; j--) {if#
\sk
\verb#(j==m-t||t>=k){Sm=0; Vt=0; c=0; for(i=1; i<=m-j; i++){if(i<k||i>=k+m-t-j)#
\sk
\verb#{Sm=Sm+S[1,i]; c++;}} if (c!=t) {sprintf(" c Error %s,%s",c,t);} Av=Sm/t;#
\sk
\verb#for(i=1; i<=m-j; i++){if(i<k||i>=k+m-t-j){Vt=Vt+(S[1,i]-Av)^2;}} sprintf#
\sk
\verb#(" tau_g = %s, T_g = [1,%s] U [%s,%s], delta_g = %s",t,k-1,k+m-t-j,#
\sk
\verb#m-j,Vt/t-(S[1,m-j]-S[1,1])/12);}} t--;}#}
\msn
This computes $\delta^{\rm Tj}_{\mhs g}$ for any {\it possible\one} $\mu$-constant deformation $g$ of $f$ whose Tjurina number is at least $\tau{-}A$. Here $A$ is a non-negative number, which can be given effectively if one knows the {\it lower bound\one} of the Tjurina number of {\it all\one} the (possible) $\mu$-constant deformations $g$ of $f$. If one is not very sure about the lower bound, it is better to set $A$ very large (for instance $\mu$), since it is replaced with a theoretically maximal number by the code. Note that the lower bound is attained on the complement of a closed analytic subset of the $\mu$-constant stratum, and one may examine the lower bound by calculating the Tjurina number at many points of the $\mu$-constant stratum. (It might be possible that the inequality \eqref{3.2} always holds in the lower bound case.) We compute the case $T_{\mhs g}\cap\If^{>\al_1+1}\eq\emptyset$ at the end if $A$ is sufficiently large. (One can see the spectral numbers by typing ``S;".) If the last outputs of all the lines are non-positive with $A$ sufficiently large, it would mean a positive answer to the conjecture. In the other case, it says nothing about the conjecture unless one can prove the existence of a $\mu$-constant deformation $g$ satisfying the desired properties. These calculations do not seem to produce a counterexample in the case the modality is at most 3 (or perhaps even 4).
\sk
It seems very difficult to prove the conjecture in the case the singularities are {\it parametrized by\one} $k\ins\N$; one can verify the conjecture by using a computer only for {\it each explicitly given\one} $k$, and never for {\it all\one} the $k$ in an algebraic way because of the ambiguity of missing spectral numbers.

\section{Newton non-degenerate case} \label{S8}
For Newton non-degenerate convenient polynomials of two variables, we can calculate the spectral numbers as in \cite{St}, \cite{Ar}, although it is not necessarily easy to determine the Tjurina spectrum (see however the code in \cite{JKSY2}). There is an exceptional case where $f$ is a linear combination of {\it three monomials.} Let
$$f\eq x_1^ax_2^b\pl x_1^c\pl x_2^d\q\h{with}\q\tfrac{a}{c}\pl\tfrac{b}{d}\slt1,\,\,\,2\less a\slt b.$$
\par\nin Set
$$\aligned&\La_0\defs\bl\{(i,j)\ins\Z^2\mid0\slt\tfrac{i}{a}=\tfrac{j}{b}\slt2\br\},\\ &\La_1\defs\bl\{(i,j)\ins\Z^2\mid0\slt\tfrac{j}{b}\slt1,\,\tfrac{i}{a}\mi\tfrac{c}{a}\slt\tfrac{j}{b}\slt\tfrac{i}{a}\br\},\\&\La_2\defs\bl\{(i,j)\ins\Z^2\mid0\slt\tfrac{i}{a}\slt1,\,\tfrac{j}{b}\mi\tfrac{d}{b}\slt\tfrac{i}{a}\slt\tfrac{j}{b}\br\},\endaligned$$
\par\nin as in a picture written in \cite{Ar}. Then the vector space $\C\{x\}/(\dd f)$ is spanned by the monomials $x^{\nu-\ob}$ for $\nu\eq(\nu_1,\nu_2)\ins\La\defs\mcup_{k=0}^2\,\La_k$ in a compatible way with the $V$-filtration, where $x^{\nu-\ob}\defs x_1^{\nu_1-1}x_2^{\nu_2-1}$ with $\ob\defs(1,1)$. This can be proved by an argument similar to the proof of the Steenbrink conjecture in \cite{JKSY1} using a resolution by a double complex of Koszul complexes as in \cite[Proposition 2.6]{Ko}, see \cite[Remark 2.1f]{JKSY1}. (It is not trivial to generalize this argument using the picture in \cite{Ar} to the case $n\gess3$, see \cite{JKSY1}.)
\sk
We see that the image of the monomial $x^{\nu-\ob}$ in $\C\{x\}/(\dd f,f)$ vanishes if $\nu\ins\La'\defs\mcup_{k=0}^2\,\La'_{k}$ with
$$\aligned&\La'_0\defs\bl\{\nu\ins\La_0\mid\nu_1\sgt a\br\},\\&\La'_1\defs\bl\{\nu\ins\La_1\mid\nu_1\sgt c\br\},\\&\La'_2\defs\bl\{\nu\ins\La_2\mid\nu_2\sgt d\br\}.\endaligned$$
\par\nin Note that $f$ mod $(\dd f)$ is represented by each of $x_1^ax_2^b$, $x_1^c$, $x_2^d$ up to a nonzero constant multiple using an argument similar to the semi-weighted-homogeneous case (since the number of monomials is three). It is easy to see that $|\La'|\eq(a{-}1)(b{-}1)$.
\sk
In the case $2b\sgt d{+}1$, we can also verify that the image of $x^{\nu-\ob}$ in $\C\{x\}/(\dd f,f)$ vanishes if $\nu\ins\La''_1$ with
$$\La''_1\defs\bl\{\nu\ins\La_1\mid\nu_1\eq c,\,\nu_2\gess d{-}b{+}1\br\}.$$
\par\nin Indeed, $x^{\nu-\ob}$ for $\nu\ins\La_1$ with $\nu_1\eq c$ can be identified up to a nonzero constant multiple with $x^{\nu+(a-c,b)-\ob}$ mod $(\dd f)$, and it belongs to $(\dd f,f)$ if $\nu_2{+}b{-}1\gess d$.
\sk
We can then conclude that the monomials $x^{\nu-\ob}$ for $\La\stm\bl(\La'\cup\La''_1)$ form a filtered basis of the filtered vector space $\bl(\C\{x\}/(\dd f,f),V\br)$ if
\begin{equation} \label{8.1}
\mu{-}\tau\eq(a{-}1)(b{-}1)\pl\max(2b\mi d\mi 1,0).
\end{equation}
Here we use a well known assertion that a filtered surjection is a filtered isomorphism if it is an isomorphism forgetting the filtration. The condition \eqref{8.1} seems to be always satisfied as far as examples are computed. Using Singular, the above argument can be implemented as follows.
\ms
\vbox{\footnotesize\sf\verb#LIB"sing.lib"; ring R=0,(x,y),ds; int a,b,c,d,i,j,p,q,t,tj; a=2;b=4;c=7;d=6;#
\sk
\verb#poly f,Av,Sm,Vt,u; u=1; if(a*d+b*c>=c*d||a>b){printf("Input Error"); exit;}#
\sk
\verb#f=x^a*y^b+x^c+y^d; tj=tjurina(f); matrix S[1][tj]; t=0; for (i=1; i<=a;#
\sk
\verb#i++) {if (i*b%a==0) {t++; S[1,t]=i*u/a;}} for (i=1; i<=c; i++) {for (j=1;#
\sk
\verb#j<b; j++) {if (b*i>a*j && (i<c || j<=d-b)) {t++; S[1,t]=i*u/c+j*(c-a)*#
\sk
\verb#u/b/c;}}} for (j=1; j<b+d; j++) {for (i=1; i<a; i++) {if (a*j>b*i && j<=d)#
\sk
\verb#{t++;S[1,t]=j*u/d+i*(d-b)*u/a/d;}}} p=1; q=1; for (i=2; i<=t; i++) {if#
\sk
\verb#(S[1,i]>S[1,p]) {p=i;} if (S[1,i]<S[1,q]) {q=i;}} if (t!=tj) {printf(#
\sk
\verb#"Serious error!");exit;} Sm=0; for(i=1; i<=t; i++){Sm=Sm+S[1,i];} Av=Sm/t;#
\sk
\verb#Vt=0; for(i=1;i<=t;i++){Vt=Vt+(S[1,i]-Av)^2;} Vt/t-(S[1,p]-S[1,q])/12;#}
\msn
Here the positive integers $a,b,c,d$ must satisfy the conditions $\tfrac{a}{c}\pl\tfrac{b}{d}\slt1$ and $2\less a\slt b$. As far as examples are computed, there are no counterexamples to the generalized Hertling conjecture among this type. However this does not imply any information about general $\mu$-constant deformations where the situation is much more complicated.
\sk
In the general Newton non-degenerate case it is not easy to determine the Tjurina spectrum (see however the code in \cite{JKSY2}). Here one has to take a monomial basis of the Jacobian ring such that it defines a section of the surjection $\C\{x\}\onto\C\{x\}/(\dd f,f)$ which is {\it compatible with the Newton filtration.} The last condition is never satisfied for arbitrary monomial bases. (This can be seen in the case where $f$ is a linear combination of three monomials with $n\eq2$ as in the above code.) It does not seem very clear whether a monomial basis given by a computer always satisfies this condition. The monomial basis given by vfilt in Singular, gmssing.lib seems to be for $\Hsf/t\Hsf$, and not for $\Hsf/\dti\!\Hsf$, for instance in the irreducible curve case with Puiseux pairs $(3,2),(1,2)$, where the $V$-filtration does not seem to be induced by ideals of $\C\{x\}$ (more precisely, by the microlocal $V$-filtration, see \cite{mic}, \cite{JKSY0}).

\section{Two Puiseux pair case} \label{S9}
In the irreducible curve case, we can determine the spectral numbers from the Puiseux pairs, see \cite{irr}. (Note that a plane curve singularity with Puiseux pairs $(k_1,n_1),(k_2,n_2)$ in the sense of \cite{irr} has Puiseux pairs $(k_1,n_1),(k_1n_2{+}k_2,n_2)$ with the definition used in Singular.) Assuming the missing spectral numbers are {\it consecutive,} one may examine whether the inequality \eqref{3.2} holds as follows.
\ms
\vbox{\footnotesize\sf\verb#LIB"sing.lib"; ring R=0,(x,y),ds; int a,b,c,d,e,i,j,k,m,p,q,r,t,mu;#
\sk
\verb#poly f,Av,Sm,Vt,u,X,Y; a=3; b=2; d=2; q=-1; r=1; c=b*q+a*r; e=a*b*d+c;#
\sk
\verb#u=1; f=(y^b-x^a)^d-x^(a*d+q)*y^r; mu=milnor(f); matrix S[1][mu]; m=0;#
\sk
\verb#for (i=1; i<e; i++) {for (j=1; j<d; j++) {X=i*u/e+j*u/d; if (X<1)#
\sk
\verb#{m++; S[1,m]=X;}}} for (i=1; i<a; i++) {for (j=1; j<b; j++) {for#
\sk
\verb#(k=0; k<d; k++) {Y=i*u/a+j*u/b; X=(Y+k)/d; if (Y<1) {m++; S[1,m]=X;}#
\sk
\verb#}}} for(i=1; i<=m; i++){p=i; for(j=i+1; j<=m; j++){if(S[1,j]<S[1,p])#
\sk
\verb#{p=j;}} X=S[1,p]; for (k=p; k>i; k--) {S[1,k]=S[1,k-1];} S[1,i]=X;}#
\sk
\verb#for (i=1; i<=m; i++) {S[1,2*m+1-i]=2-S[1,i];} m=2*m;#
\sk
\verb#if (m!=mu) {printf("Serious error!"); exit;} t=tjurina(f); Sm=0;#
\sk
\verb#for(i=1; i<=t; i++) {Sm=Sm+S[1,i];} Av=Sm/t; Vt=0;#
\sk
\verb#for(i=1; i<=t; i++) {Vt=Vt+(S[1,i]-Av)^2;} Vt/t-(S[1,t]-S[1,1])/12;#}
\msn
One may replace the positive integers $a,b,d,r$ and the integer $q$ as long as $a\sgt b\sgt r$, $ad\pl q\sgt0$, $c\defs bq\pl ar\sgt0$, and ${\rm GCD}(a,b)\eq{\rm GCD}(c,d)\eq1$, where $f$ has Puiseux pairs $(a,b),(c,d)$ in the sense of \cite{irr}. (One can see the spectral numbers by typing ``S;".) In the non-consecutive case, one may apply the last part of Remark\,\,\ref{R3.2} as in Section\,\,\ref{S7} by replacing the part after ``{\smaller\sf\verb#t=tjurina(f);#}" if it takes very long to apply the code in Section\,\,\ref{S7}.
\sk
It is rather surprising that $\mu\mi\tau$ seems to be $2$ for any positive odd number $c$ in the case $(a,b,d)\eq(3,2,2)$. Assuming $\mu\mi\tau\less 2$, we can show the inequality \eqref{3.2} rather easily in this case applying Proposition\,\,\ref{P4.1} to $T\defs\If\stm\{\mu\}$ for the case $\Cf$ is non-consecutive. Indeed, the spectral numbers are given by
$$\bl\{\tfrac{5}{12},\tfrac{11}{12},\tfrac{13}{12},\tfrac{19}{12}\br\}\cup\bl\{\tfrac{1}{2}\pl\tfrac{k}{c+12}\br\}_{k\in[1,\,c+11]}\,,$$
in particular, $\al_1\eq\tfrac{5}{12}$, $\al_2\eq\tfrac{1}{2}\pl\tfrac{1}{c+12}$, see \cite{irr}. We can then verify that $\tau\delta_T$ is expressed as
$$\aligned&2\bl(\tfrac{1}{12}\br)^2\pl\bl(\tfrac{7}{12}\br)^2\pl\tfrac{c+11}{12}\bl(1\mi\tfrac{2}{c+12}\br)\mi\tfrac{1}{c+14}\bl(\tfrac{7}{12}\br)^2\mi\tfrac{c+14}{12}\bl(\tfrac{3}{2}\mi\tfrac{1}{c+12}\mi\tfrac{5}{12})\\&=-\tfrac{c^3\pl37\one c^2\pl455\one c\pl1764}{144\one c^2\pl3744\one c\pl24192}\,,\endaligned$$
using \eqref{4.3} (where the $\al_i$ can be shifted by $-1$ to simplify the computation). Recall that the equality holds in \eqref{2.2} if $f\eq x^{c+12}$, see the last part of Section\,\,\ref{S2}. This implies the third term. We thus get the negativity of $\delta_T$, which proves the inequality \eqref{3.2} for the (possible) case $\mu\mi\tau\eq1$. One can confirm the above calculation for each odd positive integer $c$ by setting $t\eq m\mi1$ in the above code.
\sk
The inequality \eqref{3.2} in the non-consecutive case then follows from Proposition\,\,\ref{P4.1}, where the hypothesis \eqref{4.2} is shown by using the inequality (3) in \cite[Theorem~2]{JKSY2}.
\sk
In the consecutive case (where $\Cf\eq\{\mu{-1},\mu\}$), the inequality \eqref{3.2} also holds, since $\tau\delta_f^{\rm Tj}$ is expressed as
$$\aligned&2\bl(\tfrac{1}{12}\br)^2\pl\bl(\tfrac{7}{12}\br)^2\pl\tfrac{c+11}{12}\bl(1\mi\tfrac{2}{c+12}\br)\mi\bl(\tfrac{1}{2}\mi\tfrac{1}{c+12}\br)^2\mi
\tfrac{1}{c+13}\bl(\tfrac{7}{12}\pl\tfrac{1}{2}\mi\tfrac{1}{c+12}\br)^2\\&\mi\tfrac{c+13}{12}\bl(\tfrac{3}{2}\mi\tfrac{2}{c+12}\mi\tfrac{5}{12})=- \tfrac{c^4+59\one c^3+1247\one c^2+10992\one c+33840}{144\one c^3+5328\one c^2+65664\one c+269568}.\endaligned$$

\begin{rem} \label{R1.1}
It is rather interesting that no trimodal Newton-degenerate singularity seems to appear in \cite{SWZ26}; for instance the modality seems to be at least 4 if the Puiseux pairs are $(3,2),(1,3)$ with $\mu\eq42$.
\end{rem}

{\makeatletter
\def\section#1#2{}%
\def\@startsection#1#2#3#4#5#6{}%
\def\@seccntformat#1{}%
\def\refname{}%
\makeatother

\par\addvspace{1.5ex plus .5ex minus .2ex}
\noindent{\normalfont\normalsize\bfseries References}
\addcontentsline{toc}{part}{References}  
\par\vskip 1ex plus .2ex

}
\end{document}